\documentclass[12pt]{amsart}
\usepackage{amssymb,amsmath}
\numberwithin{equation}{section}
\def\T{\text}

\def\simleq{\underset\sim<}
\def\1#1{\overline{#1}}
\def\2#1{\widetilde{#1}}
\def\3#1{\widehat{#1}}
\def\4#1{\mathbb{#1}}

\def\5#1{\frak{#1}}
\def\6#1{{\mathcal{#1}}}
\def\C{{\4C}}
\def\R{{\4R}}

\def\Z{{\4Z}}

\begin{document}
\abstract
We study propagation of CR extendibility at the vertex $p$ of an analytic sector $A$ contained in a CR manifold $M$. Let $k$ be the weighted vanishing order of $M$  and $\alpha $ the complex angle of $A$ at $p$. Propagation takes place if and only if $\alpha>\frac1k$.
\endabstract
\title[Propagation of CR extendibility at the vertex]
{Propagation of CR extendibility at the vertex of a complex sector}
\author[L.~Baracco and S.~Pinton]
{Luca Baracco and Stefano Pinton }
\address{Dipartimento di Matematica, Universit\`a di Padova, via 
Trieste 63, 35121 Padova, Italy}
\email{baracco@math.unipd.it,
pinton@math.unipd.it}

\maketitle
\def\F{\mathcal F^{2,\alpha}}
\def\Giialpha{\mathcal G^{i,i\alpha}}
\def\cn{{\C^n}}
\def\cnn{{\C^{n'}}}
\def\ocn{\2{\C^n}}
\def\ocnn{\2{\C^{n'}}}
\def\const{{\rm const}}
\def\rk{{\rm rank\,}}
\def\id{{\sf id}}
\def\aut{{\sf aut}}
\def\Aut{{\sf Aut}}
\def\CR{{\rm CR}}
\def\GL{{\sf GL}}
\def\Re{{\sf Re}\,}
\def\Im{{\sf Im}\,}
\def\codim{{\rm codim}}
\def\crd{\dim_{{\rm CR}}}
\def\crc{{\rm codim_{CR}}}
\def\phi{\varphi}
\def\eps{\varepsilon}
\def\d{\partial}
\def\a{\alpha}
\def\b{\beta}
\def\g{\gamma}
\def\G{\Gamma}
\def\D{\Delta}
\def\Om{\Omega}
\def\k{\kappa}
\def\l{\lambda}
\def\L{\Lambda}
\def\z{{\bar z}}
\def\w{{\bar w}}
\def\Z{{\1Z}}
\def\t{{\tau}}
\def\th{\theta}
\emergencystretch15pt
\frenchspacing
\newtheorem{Thm}{Theorem}[section]
\newtheorem{Cor}[Thm]{Corollary}
\newtheorem{Pro}[Thm]{Proposition}
\newtheorem{Lem}[Thm]{Lemma}
\theoremstyle{definition}\newtheorem{Def}[Thm]{Definition}
\theoremstyle{remark}
\newtheorem{Rem}[Thm]{Remark}
\newtheorem{Exa}[Thm]{Example}
\newtheorem{Exs}[Thm]{Examples}
\def\Label#1{\label{#1}}
\def\bl{\begin{Lem}}
\def\el{\end{Lem}}
\def\bp{\begin{Pro}}
\def\ep{\end{Pro}}
\def\bt{\begin{Thm}}
\def\et{\end{Thm}}
\def\bc{\begin{Cor}}
\def\ec{\end{Cor}}
\def\bd{\begin{Def}}
\def\ed{\end{Def}}
\def\br{\begin{Rem}}
\def\er{\end{Rem}}
\def\be{\begin{Exa}}
\def\ee{\end{Exa}}
\def\bpf{\begin{proof}}
\def\epf{\end{proof}}
\def\ben{\begin{enumerate}}
\def\een{\end{enumerate}}
\def\dotgamma{\Gamma}
\def\dothatgamma{ {\hat\Gamma}}

\def\simto{\overset\sim\to\to}
\def\1alpha{[\frac1\alpha]}
\def\T{\text}
\def\R{{\Bbb R}}
\def\I{{\Bbb I}}
\def\C{{\Bbb C}}
\def\Z{{\Bbb Z}}
\def\Fialpha{{\mathcal F^{i,\alpha}}}
\def\Fiialpha{{\mathcal F^{i,i\alpha}}}
\def\Figamma{{\mathcal F^{i,\gamma}}}
\def\Real{\Re}
%
%
%

\section{Introduction}
Hanges and Treves prove in \cite{HT83} that a disc $A$ in a real
hypersurface $M$ of the 
complex space 
is a propagator of holomorphic extendibility across $M$. Propagation takes also
place  at a 
boundary point $p\in\partial A$ at which $\partial A$ is smooth; this is meant, though not  stated, for instance in \cite{BP85}. Thus the complex
angle for a set to be a propagator reduces 
from $2\pi$ to $\pi$. 
Going on in reduction of the angle, Zaitsev and Zampieri prove in
\cite{ZZ04} that a sector $A_\alpha$ of 
angle $\alpha\pi>\frac\pi2$ is a propagator; again, this is  not explicitely stated. 
What about $\alpha\leq\frac12$? The positive solution to the question comes from a
good balance between the 
size of $\alpha$ and the flatness of $M$ at $p$. (As for the result of \cite{ZZ04}, note
that any smooth $M$ is at least 
``$2$-flat".) For instance, consider 
the sector in $\C$
$$
A_{\frac1k+\epsilon}=\{z_1\in\C:\,\,-\pi(\frac1{2k}+\frac\epsilon2)<\T{arg}z
_1<\pi(\frac1{2k}+\frac\epsilon2)\},
$$
and the function
\begin{equation}
\label{1.1}
 h_{\frac1k+\epsilon}=
\begin{cases}
0\quad\T{ if $z_1\in A_{\frac1k+\epsilon}$},
\\
-(z_1^k+\bar z_1^k)\quad\T{ if $z_1\notin A_{\frac1k+2\epsilon}$},
\\
|z_1|^k\chi(\T{arg}(z_1))\quad\T{ if $z_1\in A_{\frac1k+2\epsilon}\setminus A_{\frac1k+\epsilon}$},
\end{cases}
\end{equation}
where $\chi\in C^\infty$ serves to connect $0$ to $-(z_1^k+\bar z_1^k)$ in the region $A_{\frac1k+2\epsilon}\setminus A_{\frac1k+\epsilon}$: it takes value $0$  (resp. $\cos(\frac\pi2+k\epsilon\pi)$) for $\T{arg}(z_1)=\pi(\frac1{2k}+\frac\epsilon2)$ (resp. $\T{arg}(z_1)=\pi(\frac1{2k}+\epsilon)$).
Let $M_{\frac1k+\epsilon}=\{z\in\C^2:\,y_2=h\}$; then  $A_\alpha\subset M_{\frac1k+\epsilon}$ for 
$\alpha={\frac1k+\epsilon}$ is a propagator of holomorphic extendibility across $M_{\frac1k+\epsilon}$ at $0$.  The proof is contained in Theorem~\ref{maintheorem} which follows but can be
much simplified in this model case. 
On the contrary, consider the smaller sector $A_{\frac1k-\epsilon}$, define $h_{\frac 1k-\epsilon}$  similarly as in \eqref{1.1} but with the pair
$A_{\frac1k+\epsilon}\subset A_{\frac1k+2\epsilon}$ replaced by 
$A_{\frac1k-\epsilon}\subset A_{\frac1k}$, and denote by $
M_{\frac1k-\epsilon}$ the hypersurface defined by 
$y_2= h_{\frac1k-\epsilon}$. It is also not restrictive to
suppose that $ h_{\frac1k-\epsilon}\geq -(z_1^k+\bar z_1^k)$.
The domain $\Omega:=\{z:\,y_2>-(z_1^k+\bar z_1^k)\}$ contains a neighborhood of the
punctured sector $A_{\frac1k-\epsilon}\setminus
\{0\}$ and also contains $\{z:\,y_2> h_{\frac1k-\epsilon}\}$.
Since $\Omega $ is pseudoconvex,  then
$A_{\frac1k-\epsilon}$ is not a propagator 
of extendibility across  $ M_{\frac1k-\epsilon}$ at $0$. This shows that the angle $\alpha\pi$ must be calibrated according to the flatness of the hypersurface.

Propagation also holds  in an exponentially degenerate boundary for an
``infinitely stretched" disc, that is, a real ray. In fact,
according to \cite{BKZ09}, the $x_1$-axis is a propagator 
in the hypersurface  $y_2=e^{-\frac1{|y_1|^s}},\,\,s\geq1$. 
This is no more true for $s<1$. On one hand the method by \cite{BKZ09} of
approximation of the real line by analytic 
discs fails, on the 
other there are superlogarithic estimates (Kohn \cite{K02}) which yield
hypoellipticity of the $\bar\partial$-Neumann 
problem. As it is classical, hypoellipticity and propagation are in contrast one to
another. This shows
that propagation in the infinitely degenerate regime is far from having
reached a complete understanding.

Instead, coming back to propagation along a nondegenerate sector on a manifold of weighted vanishing order $\geq k$, the
problem takes complete 
solution in the present paper. A sector
$A_\alpha$ for $\alpha>\frac1k$ is a propagator at the vertex.
The conclusion looks similar as to the extension at the vertex in presence of the sector property or, more generally, the rays condition by Fornaess and Rea \cite{FR85}. However, ours is propagation, and theirs, forced extension. In particular, for their result to be true, one needs to know from the beginning that the extension  holds in a neighborhood of the punctured sector whose rate at $0$ is $|z|^k$; instead, no control of the rate is required in the present paper.

Our method is largely inspired to Tumanov \cite{T94} and \cite{T97} but contains some novelties. It consists in making infinitesimal deformations of analytic discs in the Lipschitz spaces $\mathcal F^{i,\alpha}$ defined in the sequel. The crucial point is that the component of these discs which is normal to $M$ is smoothened by composition of the disc with the graphing function $h$ of $M$. When $h$ has weighted vanishing order $k$, and the tangential components are $\mathcal F^{i,\alpha}$, the normal component takes regularity $C^{k\alpha}$; in particular, if $\alpha>\frac1k$,  this is $C^1$. Another ingredient of the proof is the approximation of $\alpha$-discs by smooth discs; the convergence for the normal components is in fact in $C^{k\alpha}$. And finally, there is a smoothening argument in the implicit function theorem: a smoothening operator which gains regularity from $\mathcal F^{i,\alpha}$ to $C^{k\alpha}$ has $\mathcal F^{i,\alpha}$-inversion which is in fact in the class $C^{k\alpha}$. 

For stating our result we need to clarify the concept of weighted vanishing order. 
Let $M\subset\C^n$ be a smooth generic submanifold with $dim_{\R}(M)=2n-d$. Fix a point $p=0\in M$ and normalize coordinates so that $M$ is locally described as a graph 
$$
 M=\{ (z,w)\in\C^{d}\times\C^{n-d}\colon z=(x+iy),\, y=h(x,w)\, \} 
$$
where $h:\,\R^d\times\C^{n-d}\to\R^d$ is a smooth function such that $h(0)=0$ and $\partial h(0)=0$. Using the Taylor expansion of $h$, we can write for some $k> 2$: 
\begin{equation}
\Label{canonical}
h=P_{k-1}(w,\bar w)+xO(x^{k-1},|w|^{k-1})+O^k(w,\bar{w})
\end{equation}
where  $P_{k-1}(w,\bar{w})$ is a  homogeneous, non-harmonic polynomial of degree $k-1$.
When the graphing function $h$ has the form \eqref{canonical}, we write $h=\mathcal O^k$ and say that $M$ has weighted vanishing order $\geq k$.
This notion is coordinate free. It can be reformulated by asking, for any vector field $L\in T^{1,0}M$ $$[L^{\epsilon_1},[L^{\epsilon_2},\dots[L^{\epsilon_{k-2}},L^{\epsilon_{k-1}}]\dots]]\in T^{1,0}M\oplus T^{0,1}M$$
for any choice of $L^{\epsilon_j}$ as $L$ or $\bar{L}$.

We use the terminology ``analytic sector" of $M$ for a subset $A\subset M$ obtained as the holomorphic image under the mapping $(x+ih,w)$ of the standard sector $\{(1-\tau)^\alpha w_o:\,\,\tau\in\Delta\}$ where $\Delta$ is the standard disc and $w_o$ a vector of $\C^{n-d}$. We call $\alpha$ the angle of $A$.
Finally, we recall that a CR function $f$ is a continuous function such that $\bar L f=0$ for any $\bar L\in T^{0,1}M $.

Here is the main result of this paper.
\bt[Propagation at the vertex of an analytic sector]
\Label{maintheorem}
Let $M\subset\C^n$ be a generic submanifold of class $C^{k+3}$ with weighted vanishing order $\geq k$ at $p\in M$. We suppose that $M$ contains an analytic sector of angle $\alpha>\frac1k$ with vertex at $p$ and passing through another point $q$.

For any submanifold $M'$ of class $C^{k+3}$ with boundary $M$  at $q$, there exists a submanifold $M''\subset\C^n$ of class $C^1$ with boundary $M$  at $p$, such that any continuous CR function in $M$ that has continuous CR extension to $M'$ also has CR extension to $M''$. Moreover, given several submanifolds $M_1',\dots, M_s'$ as above in $s$ linearly indipendent directions of $T_q\C^n\setminus T_qM$, the corresponding submanifolds $M_1'',\dots,M_s''$ span $s$ linearly indipendent directions of $T_p\C^n\setminus T_pM$. 
\et 
To prove the result we first introduce in \S 2 spaces of discs with Lipschitz boundary.
Then, we proceed in three steps. The first, in \S 3, consists in showing that there exists a $\F$ disc $A'$ attached to a deformation of $M$, contained in $M\cup M'$, that has non-null normal component of the radial derivitive at $p$. The second, in \S 4, consists in showing that there exists a smooth disc $A''$, sufficiently close to $A'$ in $\F-$norm, with the radial derivative in $p$ close to those of $A'$. The third, in \S 5, consists in constructing a submanifolds over $p$ by union of the rays of a family of analytic discs obtained moving the vertices of the discs in a neighborhood of $p$.     

\section{Analytic discs with a singular boundary point} 
\subsection{The functional space $\F$}
We introduce now a particular subclass of the H\"older space $C^\alpha$ for $0 < \alpha <1$. Denote by $\tau=re^{i\theta}$ the variable in the standard disc $\Delta$. We call $\F$ the subclass of $C^\alpha$ composed by all the real continuous function $\sigma(\theta)$, $\theta\in [-\pi,\pi]$, that are $C^{2,\alpha}$ out of $0$ and for which the following norm is finite:
\begin{equation}
\Label{100}
 \| \sigma \|_{\F}=\| \sigma \|_{C^0}+\| \theta\sigma^{(1)} \|_{C^\alpha}+\| \theta^2\sigma^{(2)} \|_{C^\alpha}.
\end{equation}
(Here $\sigma^{(1)}$ and $\sigma^{(2)}$ denote first and second derivative of $\sigma$ respectively). We remark that for $\sigma\in\F$ we must have $\theta\sigma^{(1)}_{|\theta=0}=0$, for otherwise $\theta\sigma^{(1)}\to c\neq 0$, which implies $|\sigma|\geq -\frac{|c|}2\log|\theta|$ contradicting the boundedness of $\sigma$. It is easy to check that $\F$ is a Banach algebra. First we state a theorem that we use in the sequel. 

\bt[Hardy-Littlewood]
\Label{t1}
Let $\Omega$ be a bounded, lipschitz domain in $\R^n$ and let $\delta(x)$ is the distance from the boundary of $\Omega$. Let $f\in C^1(\Omega)$ and suppose that
\begin{equation}
\Label{hardy}
|\nabla f(x)|\leq C\delta(x)^{\alpha-1} 
\end{equation}
Then $f\in C^{\alpha}(\bar\Omega)$ and, moreover, the constant $C$ of \eqref{hardy} controls the $C^\alpha$-norm of $f$ by
$$ 
\|f\|_{C^\alpha}\lesssim \|f\|_{C^0}+C, 
$$ 
where ``$\simleq$" denotes inequality up to a multiplicative constant.
\et
 We denote by $\F_0$ the set of functions $\sigma$ such that $\sigma\in\F$ and $\sigma(0)=0$. Using the previous theorem, we prove an important result concerning the regularity of the composition, that we use in the sequel. Assume $(k-1)\alpha<1<k\alpha$. Let $\gamma\leq k\alpha-1$ (in particular $\gamma\leq\alpha$).  
\bp\Label{p1}
If
\begin{eqnarray}
h&\in& C^{k,\gamma}_0,\quad h\in\mathcal{O}^k,\nonumber\\
u&\in& C^{1,\gamma}_0,\quad w\in\F_0.\nonumber
\end{eqnarray}
then,
\begin{equation}
h(u,w)\in C^{1,\gamma},
\end{equation}
and
\begin{equation*}
\|h(u,w)\|_{C^{1,\gamma}}\simleq \|h\|_{C^k}\|u\|_{C^{1,\gamma}}+\|w\|_{\F}.
\end{equation*}
\ep
\bpf
We have 
$$ h(u,w)^{(1)}=h_xu^{(1)}+h_ww^{(1)}+h_{\bar w}\bar w^{(1)}. $$
Since the first term in the right is $C^{\gamma}$, then we only have to treat the second and third. Now, as for the second,
\begin{equation*}
\begin{split}
(h_ww^{(1)})^{(1)}&=h_{x,w}u^{(1)}w^{(1)}+h_{ww}(w^{(1)})^2+h_{\bar w\bar w}(\bar w^{(1)})^2+2h_{w\bar w}|w^{(1)}|^2
\\
&+h_ww^{(2)}+h_{\bar w}\bar w^{(2)}\nonumber\\
&\lesssim |\theta|^{\alpha+\alpha-1}+|\theta|^{(k-2)\alpha+2(\alpha-1)}+|\theta|^{(k-1)\alpha-2+\alpha}\nonumber\\
&\lesssim |\theta|^{-1+2\alpha}+|\theta|^{-1+\gamma}+|\theta|^{-1+\gamma}.
\end{split}
\end{equation*}
A similar estimate holds for the third term.
We then carry out our proof by the aid of the Hardy-Littlewood's theorem.

\epf

Let $T_1$ denote the Hilbert transform normalized by condition $T_1(\cdot)(1)=0$; it is easy to see that $T_1$ is a bounded operator in $\F$.
We change a little our notations and  require that for $\sigma\in\F$ we have that $\sigma(-\pi)=\sigma(\pi)$ so that $\sigma$ is naturally identified to a function of the variable $\tau=e^{i\theta}$ on the circle. In this new setting we denote by $\F_1$ the set of functions $\sigma\in\F$ such that $\sigma(1)=0$.     

\subsection{Attaching $\F$ analytic discs}
We come back to our manifold $M$ and suppose that $M$ has weighted vanishing order $\geq k$. We consider in $\C^n$ analytic discs $A(\tau)=(z(\tau),w(\tau))$, $\tau\in\Delta$, attached to $M$, that is, satisfying $A(\partial\Delta)\subset M$. If we prescribe a point $p=(z_o,w_o)$ with $y_o=h(x_o,w_o)$, and an analytic function $w_o+w_\nu(\tau)$, $\tau\in\Delta$, $\nu\in\R$, the so called CR component,  and look for an analytic completion $z_\nu(\tau)$ for $A_\nu(\tau)=p+(z_\nu(\tau),w_\nu(\tau))$ with $A_\nu(1)=p$, we are led to Bishop's equation
\begin{equation}\label{bishop_equation_1}
u_\nu(\tau)=-T_1(h(x_o+u_\nu(\tau),w_o+w_\nu(\tau)))\quad \tau\in\partial\Delta.
\end{equation}
In fact, if $u_\nu(\tau)$ solves (\ref{bishop_equation_1}), then if we set $v_\nu(\tau)=T_1u_\nu(\tau)$ and $z_\nu(\tau)=z_o+u_\nu(\tau)+iv_\nu(\tau)$, we obtain that $A_\nu(\tau)=p+(z_\nu(\tau),w_\nu(\tau))$ is holomorphic, $v_\nu(\tau)=h(x_o+u_\nu(\tau),w_o+w_\nu(\tau))$ for all $\tau\in\partial\Delta$, and finally $A_\nu(1)=p$. We consider the equation (\ref{bishop_equation_1}) in the spaces $\F$ and $C^{1,\gamma}$ for which $T_1$ is bounded.  Here we remeber two propositions that we will use subsequently.
\bp\Label{p2}
Let $h$ be of class $C^{k+2}$ and have weighted vanishing order $\geq$ $k$. Then for any $\epsilon$ there is $\delta$ such that if $\| h \|_{C^{1,\alpha}}<\delta$, $\|w_\nu\|_{\F}<\delta$, $|w_0|<\delta$, $|x_0|<\delta$ and $\nu\in\R$, then the equation (\ref{bishop_equation_1}) has a unique solution $u\in\F$ with $\|u\|_{\mathcal{F}^{2,\alpha}}<\epsilon$. The solution depends $C^k$ on $x_o$, $w_o$ and if the mapping $\nu\to W_\nu$, $\R\to\F$ is $C^k$, then it is also $C^k$ on $\nu$. Moreover, for $p=0$, $u\in C^{1,\gamma}$ where $\gamma=k\alpha-1$.
\ep
\bpf
The first claim can be proved in a standard way. In fact, consider the functional
$$F: \R^d\times\C^{n-d}\times\F\times\F\to\F,$$
defined by 
\begin{equation}
\Label{F}
F(\nu,x_0,w_0,W_\nu,u)\to u-T_1h(u+x_0,W_\nu+w_0).
\end{equation}
Then, for the partial Jacobian with respect to $u$, one has $\partial_u F\colon \dot{u}\mapsto\dot{u}-T_1\partial_xh\dot{u}$. In particular, if we evaluate at $(x_o,w_o,w_\nu,u)=(0,0,0,0)$, then this is invertible since $\partial_xh\big|_0=0$ and therefore $\partial_uF\sim\T{id}$. Thus we can apply the implicit function theorem and get the conclusion.

For the second claim, that is, $u\in C^{1,\gamma}$,  we remember that by Proposition~\ref{p1} the composition $h(u,w)$ is $C_1^{1,\gamma}$ if $u\in C_1^{1,\gamma}$, $w\in\F_1$ and $h$ is of weighted vanishing order $\geq k$ . Thus we can consider the  functional: 
\begin{eqnarray}
\Label{functional}
F^1\colon \F_1\times C_1^{1,\gamma} &\to& C_1^{1,\gamma}\nonumber\\ 
F^1(w_\nu,u)&=& u-T_1h(u, w_\nu);
\end{eqnarray} 

 This is defined as in \eqref{F} for $(x_o,w_o)=(0,0)$ but in different spaces (as an effect of the regularity of the composition).
  By the implicit function theorem, we can conclude as in the previous case that the equation $F^1=0$ has an unique solution $u\in C^{1,\gamma}_1$. On the other hand the solution of the equations $F=0$ and $F^1=0$ for the same datum $w$ are the same by uniqueness.     
\epf

\br
In the sequel we shall write $\partial_r$ for the radial derivative in the stanard disc $\Delta$.  Since by Proposition~\ref{p2} the normal component $A$ is $C^{1,\gamma}$, it makes sense to write $[\partial_r A(1)]\sim \partial_rv(1)\in (T_M \C^n)_{A(1)}$, even though $\partial_r u(1)$ may not exist.
\er

We can think of the family of discs produced by the first part of the above statement as a deformation of the disc $A(\tau)\equiv 0$ which is a trivial solution to Bishop's equation. By the next statment we show how it is possible to make infinitesimal deformations of discs which are no longer assumed to be small.

\bp
\Label{p3}
Let $h$ be as in the previous proposition, let $\tilde{w}(\tau)\in C^{1,\gamma}_1$ be small in $\mathcal{F}^{2,\alpha}$ (not necessarly in $C^{1,\gamma}$), and let $\tilde{u}(\tau)\in C^{1,\gamma}_1$ be the solution of the Bishop's equation $\tilde{u}=-T_1h(\tilde{u},\tilde{w})$. Then there exists $\delta>0$ such that for any $w_\nu(\tau)$ and $(x_o,w_o)$ with $\| w_\nu-\tilde{w} \|_{F^{2,\alpha}}\leq \delta$ and $|(x_o,w_o)|\leq\delta$ there is a unique solution $u\in C^{1,\gamma}$ of Bishop's equation 
$uT_1h(x_o+u,wo+w_\nu)$ 
with $\| u-\tilde{u} \|_{C^{1,\gamma}}<\epsilon$. Moreover, $u$ depends in a $C^k$fashin on $(xo,w_o,w_\nu)$. 
\ep
\bpf
We recall the functional 
\begin{eqnarray}
F:\R^d\times\C^{n-d}\times\F\times\F &\to& \F\nonumber\\ 
(x_o,w_o,u,w_\nu)&\to& u-T_1h(xo+u,wo+w_\nu)\nonumber
\end{eqnarray}
We know that is invertible at $(0,\tilde{w},\tilde{u})$; in particular $\partial_u F$ is injective. We can restrict $F$ to the subspace $\R^d\times\C^{n-d}\times C^{1,\gamma}\times C^{1,\gamma}$ into $C^{1,\gamma}$ and  call the new restricted functional $\tilde{F}$.  
What we have to prove is that $\partial_u\tilde{F}$ is surjective also in the restricted sense. We want to show that if $\dot{f}\in C^{1,\gamma}$ then there exists $\dot{u}\in C^{1,\gamma}$ such that $(\partial_u \tilde{F}(0,\tilde{w},\tilde{u}))\dot{u}=\dot{f}$. Since we know that $\partial_u F$ at $(0,\tilde{w},\tilde{u})$ is surjective, there exists $\dot{u}\in\F$ such that:
 $$(\partial_u \tilde{F}(0,\tilde{w},\tilde{u}))(\dot{u})= \dot{u}+T_1(h_x((\tilde{u},\tilde{w}))(\dot{u}))=\dot{f}$$      
Thus we have that:
$$\dot{u}=\dot{f}-T_1(h_x(\tilde{u},\tilde{w})(\dot{u}))$$
First we recall that $h_x(\tilde{u},\tilde{w})\in C^{1,\gamma}$. Next we observe that if $h_x(\tilde{u},\tilde{w})\in C^{1,\gamma}$, ${h_x}\big|_{0}=0$ and $\dot{u}\in\F$ then $h_x(\tilde{u},\tilde{w})\dot{u}\in C^{1,\gamma}$. Thus $\dot{u}\in C^{1,\gamma}$; this yields the invertibility of $\partial_u \tilde{F}$ at $(0,\tilde{w},\tilde{u})$ also in the space $C^{1,\gamma}$.
\epf

\section{Radial derivatives of attached $\F$ analytic discs}  
Let $A,\, M,\, M',\, p,\, q$ be as in the previous theorem. In this section we  suppose that the given disc $A$ is small in $\F$-norm (note that we don't require this hypothesis in Theorem~\ref{maintheorem}). We want to construct a family of $\F$-analytic discs $\{ A^\eta \}_{0\leq\eta\leq\epsilon}$, that depends on one real parameter $0\leq\eta\leq\epsilon$, attached to $M\cup M''$ (i.e. $\partial A^\eta\subset M\cup M''$) such that $p\in\partial A^\eta$  and  $[\partial_r A^\eta(p)]\neq 0$;
this direction depends on the direction of $M'$ normal to $M$ at $q$.

First of all, we take coordinates $(z,w)\in\C^d\times\C^{n-d}=\C^n$ with $z=x+iy$ and $w=u+iv$ in such a way that $p=0$ and

$$ U\cap M=\{ (z,w)\in\C^d\times\C^{n-d}\colon y=h(x,w) \} $$      

for $h(0)=0$ and $ \partial h(0)=0$. We also use the notation $r=y-h(x,w)$ where $r=(r_j)$, $h=(h_j)$ and $y=(y_j)$ for $j=1,\dots,d$. Let $M'$ be a manifold with boundary $M$, possibly in a neighborhood of $q$, of codimension $d-1$. This is defined, e.g., by introducing a new parameter $t\in\R^+$, and extending the domain of $h$ from $\R^d\times\C^{n-d}$ to $\R^d\times\C^{n-d}\times\R^+$ with $\partial_t h\neq 0$. Hence, $M'$ will be defined by $y=h(x,w,t)$, $t\in\R^+$.  If we consider a generic $\F$-disc $A(\tau)=(z(\tau),w(\tau))$ attached to $M\cup M'$ we call $w$ the $CR$ components of $A$ and define the $t$ components $t(\tau)$ by the equation $h(x(\tau),w(\tau),t(\tau))-y(\tau)=0$. Hence the condition $A(\tau)\in M'\setminus M$ for $\tau\in\partial\Delta$ is equivalent as to $t(\tau)>0$. We shall also let the function $t(\tau)$ depend on a small parameter $\eta\in\R^+\cup\{ 0 \}$, and denote it by $t^\eta(\tau)$. Existence of attached $\F$ discs with prescribed components $w(\tau)$ and $t^\eta$ is assured by the following statement.

\bl\label{deformation_1}
Let $h$ belong to $C^{k+2}$, $k\geq 1$ and have weighted vanishing order $\geq k$. Let $w\in\F$ and $t^\eta(\tau)\in C^{k,\alpha}$ be small. We also suppose $w(1)=0$, $t^\eta(\tau)\equiv0$ in a neighborhood of $1$ and take $(x_o,w_o)\in\R^d\times\C^{n-d}$ small. Then we can find an unique solution $u=u^\eta$ in $\F(\partial\Delta)$ of the equation

$$ u=x_o-T_1(h(u,w_o+w,t^\eta)) $$   

Moreover, if $(x_o,w_o)=(0,0)$ then $u$ is in $C^{1,\gamma}$ for $\gamma=k\alpha-1$. 
\el  
The proof of the last statement follows from the last part of Proposition~\ref{p2}.

 Using the basis $\partial r_j$, $j=1,\dots,d$, for $T_M^*\C^n$, we can identify $T_M\C^n$ (the normal bundle to $M$) to $M\times i\R^d$ by $[v]\to i(Re\left\langle \partial r_j,\, v \right\rangle)$. Let $z\in M$ and $iv^1=i\partial_t h(z)$. We have clearly $ (TM')_z=(TM)_z\oplus i\R^+v^1 $; in this case we say that $M'$ is attached to $M$ at $(z,iv^1)$ or that $M'$ is an extension of $M$ which points to the normal direction $iv^1$ at $z$.

We assume now that we are given a small $\F$ analytic disc $A$ contained in $M$ with a singularity in $p\in\partial A$ which contains another point $q$ in its boundary (here we are supposing that $p\in M$ and that $M$ has weighted vanishing order  $\geq k$). Let $q=A(-1)$ and we denote by $w(\tau)$ the CR components of $A$. Let $\partial'r$ be the square $d\times d$ Jacobian matrix of $r$ with respect to the $(z_1,\dots,z_d)$ variables. It is easy to find a real $d\times d$ matrix $G(\tau)$, $\tau\in\partial\Delta$, with $G(1)=id_{d\times d}$ and such that $G\cdot(\partial'r\circ A)$ extends holomorphically from $\partial\Delta$ to $\Delta$. To prove this, we only need to solve a linear Bishop's equation $G(\tau)=T_1(G(\tau)(\partial_x'h(u(\tau),w(\tau),t(\tau))))+id_{d\times d}$ on $\partial\Delta$ where $T_1$ is the Hilbert transform normalized by the condition $T_1(\cdot)_{|\tau=1}=0$. By means of $G$ we can define an isomorphism $(T_M\C^n)_{q}\to(T_{M}\C^n)_{p}$ which is defined, in the bases dual to $(\partial rj)_j$ by $v\mapsto G(-1)v$.

Let $\chi (\tau)$ be a real positive smooth function on $\partial\Delta$ with $\chi(-1)=1$ and whose support $supp(\chi)$ is contained in a small neighborhood of $-1$. Define $t^\eta(\tau)=\eta\chi(\tau)$ for small $\eta$ so that Lemma (\ref{deformation_1}) can be applied. Let $A^\eta$ be the family of discs of lemma (\ref{deformation_1}) with a data $w(\tau)$, $t^\eta(\tau)$ and $(x_o,wo)=(0,0)$. Let $\dot{A}$ be the derivative with respect to $\eta$ at $\eta=0$.               

\bt\Label{t5}
Let $M$ be a $C^{k+2}$ of weighted vanishing order $\geq k$, $A$ be a $\F$ in $\bar{\Delta}$, small in $\F$ norm, attached to $M$, tangent to $M$ at the singular point $p=A(1)$ and let $q$ be another point of $\partial A$, say $p=A(-1)$. Let $v^1=i(\partial_t h_j(q))_j\in i\R^d$ and $v^0=G(-1)v^1\in i\R^d$. Then 
$$ |\partial_\tau\dot{A}\big|_{1}-cv^0| < \epsilon $$   
where $c>0$ and $\epsilon$ is an error vector which can be made arbitrarly small if we correspondingly shrink $supp(\chi)$.
\et
\bpf
The proof can be found e.g. in \cite{T94}. The key point is that the CR components do not depend on $\eta$ and that the normal component of $A^\eta$ is $C^{1,\gamma}$.

\epf
Let $A^\eta=(u^\eta,w)$ be the family of $\F$ discs we defined before. We observe that in general $\partial_\tau A^\eta$ is not defined in $\tau=1$ since the CR component
 of $A^\eta$ has only regularity $\F$. However for $\eta=0$ (i.e. $A^0=A$) we can consider the function $Re\left\langle \partial r_j\circ A(\tau),\partial_\tau A^\eta(\tau) \right\rangle$ also for $\tau=1$. In fact, we have $|\partial_w r_j\circ A(\tau)|\approx |1-\tau|^{(k-1)\alpha}$ and $|\partial_\tau w(\tau)|\approx |1-\tau|^{\alpha-1}$ when $|\tau|\to 1$. This implies that $Re\left\langle \partial_w r_j\circ A(\tau),\partial_\tau w(\tau) \right\rangle\to 0$ when $|\tau|\to 1$. Thus the Taylor expansion of $\partial_\tau A^\eta$ with respect to $\eta$ gives:

$$ Re\left\langle \partial r_j\circ A,\partial_\tau A^\eta \right\rangle \big|_{1}=\eta Re\left\langle \partial r_j\circ A,\partial_\tau\dot{A} \right\rangle\big|_{1}+o(\eta)$$        

where we have used the basic hypothesis that $A$ is tangent to $M$ at $p$. Finally we observe that $\partial_\tau\dot{A}$ satisfies the conclusion of theorem (\ref{deformation_2}). It follows
$$( Re\left\langle \partial_z r_j(p),\partial_\tau A^\eta \right\rangle _{|1})_j=\eta (G_q(\partial_th_j(q))_j+\epsilon) +o(\eta).$$
where $\epsilon$ is small if $supp(\chi)$ is small. We can state the following proposition.

\bp
\Label{p10}
Let $M$ be a $C^{k+2}$ submanifold of $\C^n$ of type $\geq k$ ($k\geq 4$) and 
$A$  a small $\F$-disc attached to $M$ and tangent to $M$ at $p=A(1)$ Let $q=A(-1)$ be another point and $M'\subset\C^n$  a $C^k$-smooth submanifold with boundary $M$ at  $q$ with extra direction $v\in T_q\C^n$. Then, for any $\epsilon>0$, there exists a $C^k$-smooth family of $\F$-analytic discs $A^\eta$, $0\leq\eta\leq\eta_0$, with $A^0=A$ and $A^\eta(1)=p$  attached to $M\cup M'$ and with the property
\begin{equation}\Label{103}
[\partial_r A^\eta(1)]=\eta (G(-1)[v]+\mathcal E)+o(\eta),\quad \eta\to 0,
\end{equation} 
for some vector error $\mathcal E$.
\ep
As an immediate applications of previous proposition we obtain
\bc
Let $M\subset\C^n$ be a generic $C^{k+2}$-smooth ($k\geq 4$) submanifold through $p=0$, and let $A$ be a  small $\F$-disc  attached to $M$ and tangent to $M$ at $p=A(1)$. If $M_1',\dots,M_s'$ are $C^k$-smooth submanifolds with boundary $M$ at a point $q\in A(\partial\Delta)$ in $s$ lineary indipendent directions $[v_1],\dots,[v_s]\in (T_M\C^n)_q$, then there exist submanifolds $M_j\subset M\cup M'_j$ of class $C^k$, of the same dimension as $M$ and arbitrarly  close to $M$ in the $C^k$ norm and analytic discs $A_1,\dots,A_s$ of class $\F$ attached to $M_1,\dots,M_s$ respectively  with $A_j(1)=p$ and   such that $[\partial_r A_1(1)],\dots,[\partial_r A_s(1)]$ are linearly indipendent and are arbitrarly close to $G(-1)[v_1],\dots,G(-1)[v_s]$.   
\ec

\section{Approximation of Sectors by smooth discs}
The main result of this section is the approximation in the class $\mathcal{F}^{2,\alpha}$ by smooth discs. Let $w(\tau$ be a disc in $\mathcal F^{3,\alpha}$; define 
$$ w_\nu(\tau)= w((1-\frac{1}{\nu})\tau)-w(1-\frac{1}{\nu}); $$
note that $w_\nu$'s are smooth in $\bar{\Delta}$. Our goal is to prove
\bt\label{approximation_theorem}
We have 
$$ w_\nu(\tau)\to w \quad \textrm{in $\F(\bar{\Delta})$ for any $\alpha'<\alpha$}. $$
\et
\bpf
Let us recall that $\|\sigma\|_{\F}=\|\sigma\|_{C^{\alpha'}}+\|(1-\tau)\sigma^{(1)}\|_{C^{\alpha'}}+\|(1-\tau)^2\sigma^{(2)}\|_{C^{\alpha'}}$. Hence it suffices to show that 
\begin{equation}\Label{105}
(1-\tau)^jw_\nu^{(j)}\to (1-\tau)^jw^{(j)}\quad\textrm{in $C^{\alpha'}$;}
\end{equation} 
only for $j=0,1$ and $2$. To prove this we remark that
\begin{equation*}\Label{approximationtheorem2}
\begin{split}
 \big| ((1-\tau)^j&w^{(j)}_\nu(\tau))^{(1)} \big| \lesssim |(1-\tau)^{j-1}(w_\nu(\tau))^{(j)}|+|(1-\tau)^{j}(w_\nu(\tau))^{(j+1)}|
 \\
&\lesssim \frac{|(1-(1-\frac{1}{\nu})\tau)^{j}w^{(j)}((1-\frac{1}{\nu})\tau)|}{|(1-(1-\frac{1}{\nu})\tau)|}+\frac{|(1-(1-\frac{1}{\nu})\tau)^{j+1}w^{(j+1)}((1-\frac{1}{\nu})\tau)|}{|(1-(1-\frac{1}{\nu})\tau)|}\\
&\lesssim |(1-(1-\frac{1}{\nu})\tau)|^{\alpha-1}
\\
&\lesssim |(1-\tau)|^{\alpha-1}. 
\end{split}
\end{equation*}
Obviously, the same estimate holds with  $w_\nu$ instead of $w$, that is, we have
$$ \big| ((1-\tau)^jw^{(j)})^{(1)} \big| \lesssim C|1-\tau|^{\alpha-1}.$$
For the difference $f_\nu:=(1-\tau)^j(w_\nu^{(j)}-w^{(j)})$, we then have 
\begin{equation*}
\begin{cases}
f_\nu\to 0\quad\textrm{uniformly on compact subset of $\Delta$}\nonumber\\
|f_\nu^{(1)}|\lesssim (1-|\tau|)^{\alpha-1}.
\end{cases}
\end{equation*}
Using Hardy-Littlewood's Theorem we conclude that the sequence $f_\nu$ is bounded in $C^{\alpha}-$norm. Thus Arzela's theorem implies that $f_\nu$ converges to $0$ in $C^{\alpha'}$-norm, which is precisely \ref{105}.      
\epf

\section{Extended manifolds spanned by discs}
Let  $u_\nu$ be  $\F$-solutions  to Bishop's equations $u_\nu=-T_1h(u_\nu,w_\nu)$, and let $z_\nu=u_\nu+iv_\nu$ for $v_\nu=T_1u_\nu$. Let $u$ be the solution to $u=-T_1h(u,w)$, and set $v=T_1u$ and $z=u+iv$; suppose that $[\partial_r A(1)]=[v_o]\neq0$.
To exploit the conclusions of the preceeding section, we also assume that $w(\tau)$ has a little extra regularity, that is, $w\in \mathcal F^{3,\alpha}$; thus three derivatives are now controlled but the angle $\alpha$, which only carries geometric meaning, is unchanged.
 Since 
$$ w_\nu(\tau)\to w(\tau)\quad\textrm{in $\mathcal{F}^{2,\alpha'}(\bar{\Delta})$ for any $\alpha'<\alpha$,} $$  
and since $w_\nu(1)\equiv 0$ for all $\nu$, $z_\nu(\tau)\to z(\tau)$ in $C^{1,\beta'}(\bar{\Delta})$ by Proposition~\ref{p3} (clearly we are supposing $\alpha'$ close enough to $\alpha$ so that $\beta'\colon = k\alpha' -1>0$.) In particular for any $\epsilon>0$ and for large  $\nu$ the discs $A_\nu=(z_\nu,w_\nu)$ are in $C^{1,\beta'}$ and satisfy 
$$
|\partial_r v_\nu(1)-[v_o]|<\epsilon.
$$
 Let $\tilde{A}=(\tilde{z},\tilde{w})$ be one of these discs and recall the conclusions of Proposition~\ref{p10}. We are ready to construct a half-space $M''$ in a manifold  which contains $M$ and gains one more direction by a deformation of the disc $\tilde{A}$ such that CR functions extend from $M\cup M'$ to $M\cup M''$. For this we consider Bishop's equation $u=T_1h(x_o+u,\tilde{w}+w_o)$. According to Proposition~\ref{p3}, for any $\epsilon$ and for  small data, there is a unique solution $u$ which satisfies $\| u-\tilde{u} \|_{C^{1,\beta'}}<\epsilon$ for $\beta'<\beta\colon= k\alpha-1$. We also write $p=(x_o+ih(x_o,w_o),w_o)$, and define $A_{(x_o,w_o)}(\tau)=p+(u(\tau)+iv(\tau),\tilde{w}(\tau))$ with $v=T_1(u)$. We also write $I_{(x_0,w_0)}={A_{(x_0,w_0)}}\big|_{[-1,+1]}$ and define
 $$
 M''=\bigcup_{|x_0|\leq \delta,|w_0|\leq\delta} I_{(x_0,w_0)}([1-\epsilon,1]).
 $$
\bp\label{constructing_manifold}
$M''$ is a half space in a manifold  of codimension $d-1$ with boundary $M$ and inward conormal $v_o'$ for $v'o0$ close to $v_0$.
\ep     
\bpf
We consider the mapping
$$ \Phi\colon\C^{n-d}\times\R^d\times [1-\epsilon,1]\to V',\, (w,x,r)\to I_p(r)\quad\textrm{for $p=(x+ih(x,w),w)$.} $$
By Proposition~\ref{p2}, $\Phi$ is $C^{1,\beta'}$ in the whole of its arguments  up to $r=1$, and we have 
\begin{equation}
\Label{101}
\Phi_{(0,0,1)}'(\C^n\times\R^d\times\R^-)= T_p M+\R^+ v_0'. 
\end{equation}
In particular, $\Phi$ extends as a $C^{1,\beta'}$ mapping to a neighborhhod of $(0,0,1)$ in $\C^n\times\R^d\times(1-\epsilon,1+\epsilon)$ whose image defines a manifold $M_1=\Phi(\C^n\times\R^l\times(1-\epsilon,1+\epsilon))$. 
Thus $M''$ is a half space in this manifold; by \eqref{101} it satisfies
 $T_pM_1^+=T_pM+\R^+v_0'$.
\epf
\section{End of proof of Theorem~\ref{maintheorem}}
The last tool that we use to prove Theorem~\ref{maintheorem} is the following result.
\bt\label{t11}
Let $A$ be an analytic disc embedded in $M$ and let $p,q$ be two points in the interior. If $f$ extends at $q$ to $M'$, then it also  extends at $p$ to $M''$. Moreover, the additional directions $v'$ and $v''$ of $M'$ and $M''$ respectively are related as follows.  For a chain of small $C^1$-discs $\{A_j\}$ with $q\in\partial A_1,\dots,\, p\in\partial A_N$ and $A_j(1)=A_{j+1}(-1)$
we have 
$$ v''\sim G_N(-1)\circ\dots\circ G_2(-1)\circ G_1(-1)v' $$
\et
The proof can be found e.g. in \cite{T94}. Thus,
 propagation holds for points in the interior of $A$. By a slight modificaton of the proof it is easy to extend the result to $C^1$-points of the boundary $\partial A$. 
The result of this paper is to bring it to the singular boundary point $p=A(1)$ of an $\mathcal{F}^{2,\alpha}$ disc.

\noindent
\bpf
{\it Proof of Theorem~\ref{maintheorem}}
Let $A\in\mathcal{F}^{3,\alpha}$ with $A\subset M$. By the theorem on propagation in the interior Theorem~\ref{t11}, we can suppose that any CR continuous function in $M$ which has  CR extension to a submanifold $M'$ with boundary $M$ at $q$ also extends to a submanifold $M'''$ with boundary $M$ at $q'\in A$, where $q'$ can be choosen arbitrarly close to the vertex  and $M'$. Here, $M'$ has extra direction $v'$ and $M'''$ has extra direction $v'''$ related to $v'$ according to Proposition~\ref{p10}. (We note that given several submanifolds $M_1',\dots, M_s'$ in $s$ linearly indipendent directions in $T_q\C^n\setminus T_qM$, the corresponding submanifolds $M_1''',\dots,M_s'''$ span $s$ linearly indipendent directions in $T_{q'}\C^n\setminus T_{q'}M$.)
Let $A_1$ be a small piece of $A$ which contains $p$ as its vertex and contains $q'$ as a boundary point; moreover, $A_1$ belongs to $\mathcal F^{3,\alpha}$ and it is small in the related norm.
 Besides, we can choose  $A_1$ in such a way  that, according to the celebrated Baouendi-Treves approximation theorem,  any CR function in $M\cup M'''$ over a neighborhood of $\bar A_1$ is approximated by holomorphic polynomials. 

By Proposition~\ref{p10}, there exists an analytic disc $A_2$ attached to $M\cup M'''$ such that $[\partial_r A_2(1)]$ is close to $G(-1)[v']$. Let $w(\tau)$ be the CR component of $A_2$; we can approximate $w$ in $\|\cdot\|_{\mathcal{F}^{2,\alpha'}}$ (for any $\alpha'\leq\alpha$) by a sequence of smooth functions $w^\nu$ for which $w^\nu(1)=0$ (according to \S 3). Applying Proposition~\ref{p3}, there is a $C^{1,\beta'}$ disc $A^\nu$ attached to $M\cup M'''$ with CR component $w^\nu$ and with $\|u_2-u^{\nu_0}\|_{C^{1,\beta'}}\leq\epsilon$ (cf. \S 4). By Proposition~\ref{p3} we can construct a submanifold $M''$ with boundary $M$ at $p$ (see \S 4). Let  $f$ be a CR function over a neighborhood of $\bar A_2$ and let us  approximate $f$  by a sequence of polynomials $P_\nu$. By maximum principle there exists a subsequence $P_{\nu_\mu}$ which converges in 
$$M''=\bigcup_{|x_0|\leq \delta,|w_0|\leq\delta} I_{(x_0,w_0)}([1-\epsilon,1])$$ 
to an analytic function which is the desired extension of $f$.        
\epf

\end{document}